\documentclass[a4paper,12pt]{amsart}

\def\i{^{-1}}

\def\ci{\mathcal I}
\def\cl{\mathcal L}
\def\co{\mathcal O}

\usepackage{amssymb} 
\usepackage{latexsym} 
\usepackage{amsfonts} 
\usepackage{amsmath} 
\usepackage{eucal} 
\usepackage{bm} 
\usepackage{bbm} 
\usepackage{graphicx} 
\usepackage[english]{varioref} 
\usepackage[nice]{nicefrac} 
\usepackage[all]{xy}

\usepackage{amsthm}

\theoremstyle{plain}
\newtheorem{thm}{Theorem}[section] 
\newtheorem*{thm*}{Theorem} 
 \newtheorem{prop}[thm]{Proposition}
 \newtheorem{lem}[thm]{Lemma}
 \newtheorem{cor}[thm]{Corollary}
 \newtheorem{remark}[thm]{Remark}

\theoremstyle{definition}

\theoremstyle{remark}

\newtheorem*{remark*}{Remark}
\newtheorem*{claim*}{Claim}

\begin{document}

\title[]{Geometry of $B \times B$-orbit closures in equivariant
embeddings}
\author{Xuhua He}%
\address{School of Mathematics, Institute for Advanced Study, Princeton, NJ 08540, USA}%
\email{hugo@math.ias.edu}%
\thanks{The first author is supported by NSF grant DMS-0111298}
\author{Jesper Funch Thomsen}%
\address{Institut for matematiske fag \\ Aarhus Universitet \\ 8000 \AA rhus C, Denmark}
\email{funch@imf.au.dk}%


\begin{abstract}
Let $X$ denote an equivariant embedding of a connected
reductive group $G$ over an algebraically closed field $k$.
Let $B$ denote a Borel subgroup of $G$ and let $Z$ denote
a $B \times B$-orbit closure in $X$. When the
characteristic of $k$ is positive and $X$ is projective
we prove that $Z$ is globally $F$-regular. As a
consequence, $Z$ is normal and Cohen-Macaulay 
for arbitrary $X$ and arbitrary characteristics. 
Moreover, in characteristic zero it follows
that $Z$ has rational singularities. This
extends earlier results by the second author
and M. Brion.
\end{abstract}

\maketitle

\section{Introduction}

Let $G$ denote a connected and reductive linear algebraic group
over an algebraically closed field $k$. Let $B$ denote a Borel
subgroup of $G$. An (equivariant) embedding $X$ of $G$ is a normal
$G \times G$-variety which contains an open subset which is $G
\times G$-equivariantly isomorphic to $G$. Here we think of $G$ as
a $G \times G$-variety through left and right translation. In this
paper we study the geometry of $B \times B$-orbit closures in $X$.
Examples of such varieties include all toric varieties, all
(generalized) Schubert varieties and all large Schubert varieties
(see \cite{BP}).

The geometry of $B \times B$-orbit closures within equivariant
embeddings has been the subject of several earlier papers. 
In \cite{B} it was realized that such orbit closures were 
mostly singular with singular locus of codimension 2. 
In the special case of the {\it wonderful compactification} 
of a semisimple group $G$ of adjoint type, this was later 
strengthened in \cite{BP}, where it was proved that 
closures of orbits of the form $B g B$, for $g \in G$, are normal and
Cohen-Macaulay. Closures of this form are called {\it large
Schubert varieties}. Using the concept of global 
$F$-regularity the latter result was generalized to arbitrary 
$X$ and $G$ in \cite{BT}. For arbitrary 
$B \times B$-orbit closures it seems that normality
and Cohen-Macaulayness is only known for the
wonderful compactifications \cite[Rem.1]{B2}. 
In the present paper we show that all $B
\times B$-orbit closures for  arbitrary  $X$ and $G$
will be normal and Cohen-Macaulay.  Moreover, when 
the field $k$ has characteristic $0$ we will show 
that such orbit closures have rational singularities.
As in \cite{BT} the main technical tool will be that 
of  global $F$-regularity.

Global $F$-regularity was introduced by
K. Smith in \cite{S}. By definition a
projective variety $Z$ over a field
of positive characteristic is globally
$F$-regular if every ideal of some
homogeneous coordinate ring of $Z$
is tightly closed. Any globally
$F$-regular variety will be normal
and Cohen-Macaulay. Moreover, every
homogeneous coordinate ring of $Z$
will share the same properties.
Another consequence is that the higher
cohomology groups of nef line bundles
on $Z$ will be zero. Known
classes of globally $F$-regular
varieties include projective toric varieties
\cite{S}, (generalized) Schubert varieties
\cite{LPT} and projective large Schubert
varieties \cite{BT}. In this paper we prove
that every $B \times B$-orbit closure in a
projective embedding $X$ of a reductive
group $G$ is globally $F$-regular. Notice
that varieties of this form include the
mentioned classes above.

The paper is organized as follows. In Section 2 we introduce
notation. In Section 3 we give a short introduction to Frobenius
splitting, canonical Frobenius splitting and global
$F$-regularity. In section 4 we present the main technical result
(Proposition \ref{mainlemma}) which relates the mentioned concepts
from Section 3. In section 5 we describe the $G \times G$-orbit
closures in a toroidal embedding. Section 6 describes the
decomposition of the closure of a $B \times B$-orbit into the
union of some $B \times B$-orbits for toroidal embeddings. This is
a generalization of Springer's result in \cite{Sp} on the
wonderful compactification. As a by-product of this description we
obtain, that any Frobenius splitting of a toroidal embedding $X$
which compatibly Frobenius splits the boundary components and the
large Schubert varieties of codimension 1, will automatically
compatibly Frobenius split all $B \times B$-orbit closures in $X$.
This is used in Section 7 to conclude that all $B \times B$-orbit
closures in a toroidal embedding are simultaneous canonical
Frobenius split. In section 8 we prove that any $B \times B$-orbit
closure in a projective embedding (over a field of positive
characteristic) is globally $F$-regular. The proof of this
proceeds by reducing to the case when $X$ is toroidal and then
using the results of the previous sections. Finally in Section 9
we treat the characteristic $0$ case by descending the results 
from Section 8 to positive characteristic.

\section{Notation}

Throughout this paper $G$ will denote a connected reductive
linear algebraic group over an algebraically closed field
$k$. The associated semisimple and connected group of 
adjoint type will be denoted by $G_{\rm ad}$. The 
associated canonical morphism is denoted by $\pi_{\rm ad} : 
G \rightarrow G_{\rm ad}$. We will fix a maximal torus 
$T$ and a Borel subgroup $B \supset T$ of $G$. 

The set of roots determined by $T$ will be denoted by $R$ and we
define the subset of positive roots $R^+$ of $R$ to be the set of
roots $\alpha \in R$ such that the $\alpha$-weight space of the
Lie algebra of $B$ is nonzero. The set $\Delta = \{ \alpha_1,
\dots, \alpha_l \} $ of simple positive roots will be indexed by
$I = \{ 1, \dots, l \}$. For each subset $J \subset I$ we let $P_J
\supset B$ denote the corresponding parabolic subgroup of $G$. The
associated Levi subgroup containing $T$ will be denoted by $L_J$
while we use the notation $U_J$ to denote the unipotent radical of
$P_J$. The notation $U_J^-$ and $L_J^-$ will be used for the 
equivalent subgroups in the parabolic subgroup $P_J^-$ opposite 
to $P_J$. When $J$ is empty we simple denote $P_J^-$ by $B^-$ 
and $U_J$ by $U$. The semisimple group of adjoint type associated 
with $L_J$ is denoted by $G_J$.

To each root $\alpha \in R$ there is an associated reflection
$s_\alpha$ in the Weyl group $W = \nicefrac{N_G(T)}{T}$. The
reflection associated with the simple root $\alpha_i$ is called
simple and will be simply written as $s_i$. We may then write each
element $w$ in $W$ as a product of simple reflection and the
minimal number of factors in such a product is the length of $w$
and will be denoted by $l(w)$. The unique element of maximal
length will be denoted by $w_0$. For $J \subset I$, we denote by
$W_J$ the subgroup of the Weyl group $W$ generated by the simple
reflections $s_i$ for, $i \in J$, and by $W^J$ the set of minimal
length coset representatives of $W/W_J$. The element in $W_J$ of
longest length is denoted by $w_0^J$. For an element $w \in W$ we
let $\dot{w}$ denote a representative for $w$ in the normalizer
of $T$. Moreover, we define $R(w) = \{ \alpha \in R^+ : w  \alpha \in R^+ \}, 
$ and denote by $U_w$ the subgroup of $B$ generated by the root
subgroups $U_\alpha$ for $\alpha \in R(w)$. Then we let $B_w$
denote the subgroup $T U_w$ of $B$.

By a variety over $k$ we mean a reduced and
separated scheme of finite type over $k$. In
particular, a variety need not be irreducible.

\section{Generalities on Frobenius splitting}
\label{F-splitting}

Let $X$ be a scheme of finite type over an algebraically
closed field $k$ of positive characteristic $p>0$.
The {\it absolute Frobenius morphism} $F : X \rightarrow
X$ on $X$ is the morphism of schemes which on the level
of points is the identity map and where the associated map of
sheaves
$$F^\sharp : \mathcal O_{X} \rightarrow F_* \mathcal O_{X},$$
is the $p$-th power map. A {\it Frobenius splitting} of $X$ is an
$\mathcal O_{X}$-linear morphism
$$ s : F_* \mathcal O_{X} \rightarrow \mathcal O_{X},$$
such that the composition $s \circ F^\sharp$ is the identity map.

\subsection{Compatibly split subschemes}

Let $Y$ denote a closed subscheme of $X$ with sheaf
of ideals $\mathcal I_Y$. A Frobenius splitting $s$
of $X$ is said to {\it compatibly Frobenius split} $Y$
if $s(\mathcal I_Y) \subset \mathcal I_Y$. In this
case there exists an induced Frobenius splitting
of $Y$. When $Y$ is compatibly Frobenius split
by $s$ then any irreducible component of $Y$ will
also be compatibly Frobenius split by $s$. Moreover,
if $Y'$ is another (by $s$) compatibly Frobenius split
closed subscheme then the scheme theoretic intersection
$Y \cap Y'$ will also be compatibly Frobenius split
by $s$.

\subsection{Push-forward}

Let $f : X \rightarrow X'$ denote a morphism of
schemes of finite type over $k$. Assume that $X$
admits a Frobenius splitting $s$ which compatibly
splits a closed subscheme $Y$. If the induced
map $f^\sharp :  \mathcal O_{X'} \rightarrow
f_* \mathcal O_{X},$ is an isomorphism, then
$s$ induces by push-forward a Frobenius splitting
of $X'$ which compatibly Frobenius splits the 
scheme theoretic image of $Y$.

\subsection{Stable Frobenius splitting along divisors}

Let $D$ denote an effective Cartier divisor on $X$ and
let $s_D$ denote the canonical section of the associated
line bundle $\mathcal O_X(D)$. Then $X$ is said to admit
a stable Frobenius splitting along $D$ if there exists 
a positive integer $e$ and an $\mathcal O_{X}$-linear morphism
$$ s : F_*^ e \mathcal O_{X}(D)  \rightarrow \mathcal O_{X},$$
such that $s(s_D)=1$. Notice that in this case the composition
of $s$ with the canonical map
$ \mathcal O_{X} \rightarrow F_*^e \mathcal O_{X}(D)$,
defined by $s_D$,  is a Frobenius splitting of $X$.
If $D'$ is another effective divisor then it
is known (see e.g. \cite[Lemma 3.1]{BT}) that $X$ is stably Frobenius split
along the sum $D+D'$ if and only if $X$ is stably
Frobenius split along both $D$ and $D'$.

When $X$ admits a stable Frobenius splitting $s$
along $D$ and $Y$ is a closed subscheme of $X$,
then we say that $s$ compatibly Frobenius
splits $Y$ if
$ s\big(  F_*^ e \big( \mathcal I_Y \otimes \mathcal O_{X}(D)
\big) \big ) \subset  \mathcal I_Y$ and, moreover,
none of the components of $Y$ are contained in
the support of $D$.

\subsection{Canonical Frobenius splitting}

Let now $G$ be a connected reductive linear
algebraic group. Fix a Borel subgroup $B$
and a maximal torus $T \subset B$ of $G$.
When $X$ is a $B$-variety there is an induced
action of $B$ on the set of $\mathcal O_X$-linear
maps ${\rm Hom}_{\mathcal O_X} (F_* \mathcal O_{X},
\mathcal O_{X})$. More precisely, when $b \in B$
and $f \in \mathcal O_X(V)$, for $V$ open in $X$,
then we define $b \cdot f$ to be the function on
$b V$ defined by $(b \cdot f)(v) = f(b^{-1} v)$.
Then for $s : F_* \mathcal O_{X} \rightarrow
\mathcal O_{X}$ we define $(b \star s) : F_* \mathcal O_{X}
\rightarrow  \mathcal O_{X}$ by
$$(b \star s) (f) = b \cdot s( b^{-1} \cdot f).$$
We regard ${\rm Hom}_{\mathcal O_X}
(F_* \mathcal O_{X}, \mathcal O_{X})$ as a
$k$-vectorspace by letting $z \in k$ act
on $s : F_* \mathcal O_{X} \rightarrow
\mathcal O_{X}$ as
$$ (z . s) (f) = z^{\nicefrac{1}{p}}
s(f).$$
We may then define the following important
concept : a Frobenius splitting $s$ of $X$
is said to be $(B,T)$-canonical if :
\begin{itemize}
\item $t \star s = s$,  $\forall t \in T$.
\item Let $\alpha \in \Delta$ and let
$x_\alpha : k \rightarrow G$ be the associated
homomorphism of algebraic groups satisfying $t
x_\alpha(z) t^{-1} =  x_\alpha(\alpha(t) z)$, $t \in T$.
Then
$$ x_\alpha(z) \star s = \sum_{i=1}^{p-1} z^i .
s_{i, \alpha} ~~, \text{for all } z \in k,$$ for certain
fixed $s_{i, \alpha} \in {\rm Hom}_{\mathcal O_X} (F_* \mathcal O_{X}, \mathcal O_{X})$.
\end{itemize}

When $X$ is a $B$-variety we define the variety $G \times_B X$
to be the quotient of $G \times X$ by the $B$-action defined
by $b . (g,x) = (g b^{-1}, b x)$ for $b \in B, g \in G$ and
$x \in X$. With this notation we have the following crucial
result connected with canonical Frobenius splittings (see
e.g. \cite[4.1.E(4)]{BK})

\begin{prop}
\label{canonical}
Let $X$ be a variety admitting a $(B,T)$-canonical Frobenius
splitting $s$. Then the variety $G \times_B X$ admits a
$(B,T)$-canonical Frobenius splitting such that
$ \overline{B\dot{w}B} \times_B X$ is compatibly Frobenius split for all
$w \in W$ and such that $G \times_B Y$ is compatibly
Frobenius split
for all $B$-stable subvarieties of $X$ which are compatibly
Frobenius split by $s$.
\end{prop}

\subsection{Strong $F$-regularity} A general reference
for this subsection is \cite{HH}.
Let $K$ be a field of positive characteristic $p>0$ and 
let $R$ denote a commutative $K$-algebra essentially of
finite type, i.e. equal to some localization of a finitely 
generated $K$-algebra. We say that $R$ is {\it strongly $F$-regular} 
if for each $s \in R$, not contained in a minimal prime of $R$, 
there exists a positive integer $e$ such that the $R$-linear 
map $F_s^e: R \rightarrow F_*^e R$,
$r \mapsto r^{p^e} s$, is split. When $R$ is strongly 
$F$-regular then $R$ is normal and Cohen-Macaulay. 
Moreover, all ideals in $R$ will be tightly closed and 
thus $R$ will be $F$-rational, i.e. every parameter 
ideal is tightly closed.

The ring $R$ is strongly $F$-regular if and only if all of 
its localized rings are strongly $F$-regular. Thus, we define 
a scheme $X$ of finite type over $K$ to be strongly $F$-regular 
if all of its local rings $\mathcal O_{X,x}$, for $x \in X$,
are strongly $F$-regular. Then the affine scheme ${\rm Spec}(R)$ 
(when $R$ is a finitely generated $K$-algebra) is strongly 
$F$-regular precisely when $R$ is strongly $F$-regular. 

\subsection{Global $F$-regularity}

Consider an irreducible projective variety $X$ over $k$.
For an ample line bundle $\cl$ on $X$ we define the
associated section ring to be
$$R = R(X, \cl) := \bigoplus_{n \in \mathbb{Z}} \Gamma(X, \cl^n).$$
We then say that $X$ is {\it globally $F$-regular} if
the ring $R(X, \cl)$ is strongly $F$-regular for some (or
equivalently, any) ample invertible sheaf $\cl$ on $X$.
Global $F$-regularity was introduced by K. Smith in 
\cite{S}. When $X$ is globally $F$-regular then $X$ is 
also strongly $F$-regular. In particular, globally $F$-regular 
varieties are normal, Cohen-Macaulay and locally $F$-rational.

The following important result by Smith \cite[Theorem 3.10]{S}
connects global $F$-regularity, Frobenius splitting and strong
$F$-regularity.

\begin{thm}
\label{Thm-Smith} If $X$ is an irreducible
projective variety over $k$ then the following
are equivalent:

\begin{enumerate}
\item $X$ is globally $F$-regular.

\item $X$ is stably Frobenius split along an ample effective
Cartier divisor $D$ and the (affine) complement $X \setminus D$ is
strongly $F$-regular.

\item $X$ is stably Frobenius split along every effective Cartier
divisor.
\end{enumerate}
\end{thm}

The connection between (1) and (3) in this theorem leads to the
following result which can be found in \cite{LPT}.

\begin{cor}
\label{pushforward} Let $f : \tilde{X} \rightarrow X$ be a
morphism of projective varieties. If the connected map $f^\sharp : \co_X
\rightarrow f_* \co_{\tilde{X}}$ is an isomorphism and $\tilde{X}$
is globally $F$-regular then $X$ is also globally $F$-regular.
\end{cor}

\section{Some criteria for globally F-regularity}

Throughout this section we assume that $k$ has positive
characteristic. The following result connects canonical
Frobenius splitting and global $F$-regularity.

\begin{prop}
\label{mainlemma}
Let $Y$ be an irreducible projective $B$-variety. Let
$y \in Y$ and $w \in W$. Define $Y' = Y - B \cdot y$ and
assume that
\begin{enumerate}
\item $B_w \cdot y =B \cdot y$ and $B \cdot y$ is dense in $Y$.
\item $Y$ admits a
$(B,T)$-canonical Frobenius splitting which compatibly splits the
subvariety $Y'$.
\item $Y$
is strongly $F$-regular.
\end{enumerate}
Write $w= s_{i_1} s_{i_2} \cdots s_{i_n}$ as a reduced
product of simple reflections and define
$$ Z = P_1 \times_B P_2 \times_B \cdots
\times_B P_n \times_B Y,$$ where $P_j = B \cup B \dot s_{i_j} B$ is a
minimal parabolic subgroup. Then $Z$ is globally $F$-regular.
\end{prop}
\begin{proof}
Let $\mathcal L$ denote an ample line bundle on $Z$. Since $Y$ is strongly
$F$-regular, $Y$ is normal. Moreover the Picard
group of $B$ is trivial. Thus we may consider $\cl$ as a
$B$-linearized line bundle. In particular, $B$ acts linearly on
the finite dimensional vector space ${\rm H}^0(Z, \mathcal L)$ of
global sections of $\cl$ and we may thus 
find a nonzero global section $s$ which is $B$-invariant 
up to scalars.

Let $z = [\dot s_{i_1}, \dots, \dot s_{i_n}, y] \in Z$. Then by
assumption (1) the orbit $B \cdot z$ is dense in $Z$ with
complement equal to the union of the subsets
$$ Z_i = P_1 \times_B  \dots \times_B B \times_B
\dots \times_B P_n \times_B Y, ~i=1, \dots, n,$$
$$ Z'_j =  P_1 \times_B P_2 \times_B \cdots
\times_B P_n \times_B Y'_j,~ j=1, \dots,m ,$$ where $Z_i$ is
defined by substituting $B$ with $P_i$ in the definition of $Z$
and $Y'_j$, $j=1, \dots,m$, denotes the components of $Y'$.
As the support ${\rm supp}(s)$ of $s$ is $B$-stable and of
codimension 1 in $Z$ it follows that  ${\rm supp}(s)$ is
contained in  $Z - B \cdot z$, i.e in
the union of $Z_i$, $i=1, \dots,n$ and $Z'_j$, $j=1, \dots, m$. In
particular, we may choose nonnegative integers $n_i$ and $m_j$
such that the zero divisor of $s$ in $Z$ equals
$$Z(s) = \sum_{i=1}^n  n_i Z_i +   \sum_{j=1}^m m_j Z'_j.$$

By assumption (2) and Proposition \ref{canonical}
the variety $Z$ admits a Frobenius splitting
which compatibly Frobenius splits $Z'_j$, $j=1, \dots, m$
and $Z_i$, $i=1, \dots, n$. Let $Y^0$ denote the ($B$-invariant)
nonsingular locus in $Y$. As $Y$ is normal the complement $Y -
Y^0$ is of codimension $\geq 2$. Now define
$$ Z^0 =  P_1 \times_B P_2 \times_B \cdots
\times_B P_n \times_B Y^0.$$ Then $Z^0$ is a smooth variety which
admits a Frobenius splitting compatibly splitting the divisors
$  Z_i \cap Z^0$, $i=1, \dots, n$ and the subvarieties
$Z'_j \cap Z^0$, $j =1, \dots, m$.
As $Z^0$ is smooth this implies (see e.g. \cite[Lemma 1.1]{LPT})
that $Z^0$ admits a stable Frobenius splitting along the effective
Cartier divisor :
$$ \sum_{i=1}^n  (Z_i \cap Z^0) +   \sum_{j=1}^m \delta_j (Z'_j
\cap Z^0),$$
where $\delta_j=0$ if $Z'_j$ is not a divisor and else $\delta_j=1$.
As a consequence, $Z^0$ admits a  stable Frobenius splitting along
$$ \sum_{i=1}^n  n_i(Z_i \cap Z^0) +   \sum_{j=1}^m m_j(Z'_j
\cap Z^0).$$

In other words, $Z^0$ is stably Frobenius split along
the Cartier divisor defined by the restriction of $s$ to $Z^0$.
Thus the morphism
$$ \mathcal O_{Z^0} \rightarrow F_*^e
\mathcal O_{Z^0}(Z(s) \cap Z^0), $$
defined by the restriction of $s$ to $Z^0$ splits for some
sufficiently large integer $e$.

As $Y$ is normal so is $Z$. Moreover, $Z - Z^0$ has codimension
$\geq 2$ and thus $i_* i^* \mathcal M$ for any line bundle
$\mathcal M$ on $Z$ where $i$ denotes the inclusion map of $Z^0$
in $Z$. Applying the functor $i_*$ to the stable splitting above
we find that $Z$ admits a stable Frobenius splitting along the
effective Cartier divisor defined by $s$. Moreover, as $Y$ is
strongly $F$-regular also $Z$ and hence $Z - {\rm supp}(s)$ is
strongly $F$-regular (see e.g. \cite[Lemma 4.1]{LS}).
This proves that $Z$ is globally $F$-regular
and ends the proof.
\end{proof}

For convenience of the reader we include the following
result (see \cite[Lemma 2.11]{R}) which we will use
in the proof of the next proposition.

\begin{lem}
\label{Kempf}
Let $f : X \rightarrow Y$ denote a projective morphism
of irreducible varieties and let $X'$ denote a closed
irreducible subvariety of $X$. Consider the image
$Y' = f(X')$ as a closed subvariety of $Y$. Let
$\mathcal L$ denote an ample line bundle on $Y$
and assume
\begin{enumerate}
\item $f_* \mathcal O_X = \mathcal O_Y$.
\item ${\rm H}^i(X, f^* \mathcal L^n ) =
 {\rm H}^i(X', f^* \mathcal L^n )  = 0$
for $i>0$ and $n>>0$.
\item The restriction map
${\rm H}^0(X, f^* \mathcal L^n ) \rightarrow
{\rm H}^0(X', f^* \mathcal L^n )  = 0$ is surjective
for $n>>0$.
\end{enumerate}
Then the induced map $f' : X' \rightarrow Y'$
is a rational morphism, i.e. $f'_* \mathcal O_{X'} =
\mathcal O_{Y'}$ and ${\rm R}^i f'_* \mathcal
O_{X'} = 0$, $i>0$.
\end{lem}

\begin{prop}
\label{rational}
Let $X$ denote an irreducible $G$-variety and let $Y$ denote a
closed irreducible $B$-subvariety of $X$. Assume that $X$
admits a $(B,T)$-canonical Frobenius splitting which compatibly
splits $Y$. Let $P_1, \dots, P_n$ denote a collection of
minimal parabolic subgroups of $G$. Then the natural
map
$$ f : Z = P_1 \times_B \cdots \times_B P_n \times_B
Y \rightarrow (P_1 \cdots P_n) \cdot Y \subset X,$$
is a rational morphism, i.e. $R^if_* \mathcal O_{Z} = 0$, $i>0$,
$f_* \mathcal O_Z = \mathcal O_{f(Z)}$.
\end{prop}
\begin{proof}
Define $Z_X =  P_1 \times_B \cdots \times_B P_n \times_B X$.
As $X$ is a $G$-variety we may identify $Z_X$ with the product
$Z(P_1,\dots,P_n) \times X$, where
$Z(P_1,\dots,P_n)$ denotes the Bott-Samelson variety
$P_1 \times_B \cdots \times_B P_n /B$. We 
define $g :Z_X \rightarrow X$ to be the associated
projection map.  As $Z(P_1,\dots,P_n)$ is an irreducible
projective variety we have $g_* \mathcal O_{Z_X} =
\mathcal O_X$.

Let $Z_{X,i}$, $i=1, \dots, n$, denote the Cartier divisor
$$Z_{X,i} =  P_1 \times_B \cdots \times_B B \times_B
\cdots \times_B P_n \times_B X, $$
in $Z_X$, where $P_i$ in the definition of $Z_X$ is
substituted by $B$. Then, by Proposition \ref{canonical}, the
variety $Z_X$ admits a Frobenius splitting $s$
which compatibly splits the subvariety $Z$ and the divisors
$Z_{X,i}$, $i=1, \dots, n$. Thus by \cite[Lem.1.1]{LPT} the
Frobenius splitting $s : F_* \mathcal O_{Z_X} \rightarrow
O_{Z_X}$ maps through the morphism
$$  F_* \mathcal O_{Z_X} \rightarrow
F_* \Big( \mathcal O_{Z_X}\Big(\sum_{i=1}^n Z_{X,i}\Big ) \Big)$$
defined by the product of the canonical sections
of the Cartier divisors $Z_{X,i}$, $i=1, \dots, n$.
Thus we may regard $s$ as a stable Frobenius splitting
of $X$ along $\sum_{i=1}^n Z_{X,i}$ which compatibly
splits $Z$. By \cite[Lem.4.3, Lem.4.4]{T} we conclude
that $Z_X$ admits a stable Frobenius splitting along
any divisor of the form
$$ \sum_{i=1}^n n_i Z_{X,i},$$
with $n_i$ being positive integers, which compatibly
Frobenius splits $Z$.

Let $\mathcal L$ denote any ample line bundle on $X$.
Choose $n_i$, $i=1, \dots,n$, such that the line bundle
$$\mathcal L'_m= g^*\mathcal L^{p^m} \otimes
\mathcal O_{Z_X}\big(\sum_{i=1}^n n_i Z_{X,i}\big)$$ is ample on
$Z_X$ for all $m>0$ (that this is possible follows e.g. from
\cite[Lem.6.1]{LT}). By \cite[Lem.4.8]{T} there exists, for some
$m$, an embedding of abelian groups
$$ {\rm H}^j\big(Z_X, \mathcal I_Z \otimes g^*\mathcal L \big)
\subseteq  {\rm H}^j\big(Z_X, \mathcal I_Z \otimes \mathcal L'_m \big),$$
for all $j$. So by \cite[Thm.1.2.8]{BK} and the ampleness of
$\mathcal L_m'$ it follows that $ {\rm H}^j\big
(Z_X, \mathcal I_Z \otimes g^* \mathcal L \big)$ is zero for $j>0$.
Similarly (with $Z$ substituted with $Z_X$) we may
conclude that ${\rm H}^j\big (Z_X, g^* \mathcal L \big)$ is zero
for $j>0$. Together these two latter statements imply
that ${\rm H}^j\big (Z, g^* \mathcal L \big)$ is also zero
for $j>0$.

Applying Lemma \ref{Kempf} now ends the proof.
\end{proof}

Combining the two propositions above with Corollary
\ref{pushforward} we find

\begin{thm}
\label{glregthm}
Let $X$ denote an irreducible $G$-variety and let $Y$ be a closed
irreducible $B$-subvariety of $X$. Assume that $X$ admits a
$(B,T)$-canonical Frobenius splitting which compatibly splits $Y$.
Let $y \in Y$ and $w \in W$ and assume that the triple $(Y,y,w)$
satisfies the assumptions in Proposition \ref{mainlemma}. Then
$(\overline{B\dot{w}B})Y$ is globally $F$-regular.
\end{thm}

\section{The $G \times G$-orbit closures in toroidal embeddings}

Consider $G$ as a $G \times G$-variety by left and right
translation. An {\it equivariant
$G$-embedding} (or simply a $G$-embedding) is a normal
$G \times G$-variety $X$ containing an open subset
which is $G \times G$-equivariantly isomorphic to $G$.

\subsection{Wonderful compactifications}
\label{emb}

When $G=G_{\rm ad}$ is of adjoint type there exists a
distinguished equivariant embedding ${\bf X}$ of
$G$ which is called the {\it wonderful compactification}
(see e.g. \cite[6.1]{BK}).

The boundary ${\bm X} - G$ of ${\bm X}$ is a union of irreducible
divisors ${\bm X}_i$, $i \in I$, which intersect transversally. For
a subset $J \subset I$ we denote the intersection $\cap_{j \in J}
{\bm X}_j$ by ${\bm X}_J$. Then ${\bm Y} := {\bm X}_I$ is the
unique closed $G \times G$-orbit in ${\bm X}$. As a $G \times
G$-variety ${\bm Y}$ is isomorphic to $\nicefrac{G}{B} \times
\nicefrac{G}{B}$.

\subsection{Toroidal embeddings}
\label{toi-emb}

An embedding $X$ of a reductive group $G$ is called {\it toroidal}
if the canonical map $\pi_{\rm ad} : G \rightarrow G_{\rm ad}$ admits an
extension $\pi : X \rightarrow  {\bf X}$ into the wonderful
compactification ${\bf X}$ of the group  $G_{\rm ad}$ of adjoint
type.

\subsection{The $G \times G$-orbits} For the rest of this section
we assume that $X$ is a toroidal embedding of $G$. The boundary
$X-G$ is of pure codimension 1 (see \cite[Prop.3.1]{Hartshorne}). 
Let $X_1, \dots, X_n$ denote the boundary divisors. For 
each $G \times G$-orbit closure $Y$ in $X$ we then associate 
the set 
$$K_Y = \{ i \in \{ 1, \dots, n \} \mid Y \subset X_i \},$$ 
where by definition $K_Y = \varnothing$ when $Y = X$.
Then by \cite[Prop.6.2.3]{BK}, $Y = \cap_{i \in K_Y} X_i$. 
Moreover, we define 
$$\ci = \{ K_Y \subset \{ 1, \dots, n \} \mid Y \text{~a
$G \times G$-orbit closure in $X$ } \},$$
and write $X_K := \cap_{i \in K} X_i$ for $K \in \ci$.
Then $(X_K)_{K \in \ci}$ are the closures of $G \times G$-orbits
in $X$. When $X$ is the wonderful compactification of $G_{\rm ad}$
then $\ci=\mathcal P(I)$,  where $\mathcal P(I)$ denotes the
set of subsets of $I$. Moreover, the $G \times G$-equivariant
map $\pi: X \rightarrow \bm X$ induces a map
$p: \ci \rightarrow \mathcal P(I)$ such that $\pi(X_K)={\bm
X}_{p(K)}$.

\subsection{The base points}
\label{basepoints}
Let $X'$ denote the closure of $T$ within $X$ and let similarly
${\bf X}'$ denote the closure of $T_{\rm ad}= \pi_{\rm ad}(T)$ 
within ${\bf X}$.
Let $X_0$ denote the complement of the union of the closures
$\overline{B \dot{s}_iB^-}$, $i=1, \dots,l,$ within $X$.
Then $X_0$ is an open $B \times B^-$-stable subset of $X$.
Moreover, if we let $X_0'$ denote the
intersection of $X'$ and $X_0$ then the map
$$ U \times U^- \times X'_0 \rightarrow X_0,$$
$$ (u,v,x) \mapsto (u,v) z.$$
is an isomorphism (see  \cite[Prop.6.2.3(i)]{BK}). With similar
definitions for $\bf X$ we also obtain an isomorphism
$$ U \times U^- \times {\bf X}'_0 \rightarrow {\bf X}_0.$$
The above defined subsets are related in the way that
$\pi^{-1}({\bf X}'_0) =  X'_0$ and consequently also
$\pi^{-1}({\bf X}_0) =  X_0$.

The set ${\bf X}_0'$ is a toric variety (with respect to $T_{\rm
ad}$). In particular, it contains finitely many $T \times
T$-orbits. The $T \times T$-orbits are classified by the set 
$ \mathcal P(I)$ of subsets of $I$. We may choose representatives
$\bm h_J$, $J \subset I$, for these orbits such that $\bm h_J$ is
invariant under the groups $U_{I-J}^- \times U_{I-J}$, ${\rm
diag}(L_{I-J})$ and $Z(L_{I-J}) \times Z(L_{I-J})$ (see e.g.
\cite[1.1]{Sp}). Such a representative $\bm h_J$ 
is then uniquely determined.

Each $G \times G$-orbit in ${\bf X}$ intersects ${\bf X}_0'$ in a
unique $T \times T$-orbit (see \cite[Prop.6.2.3(ii)]{BK}). In
particular, the elements $\bm h_J$ are also representatives for the $G
\times G$-orbits in ${\bf X}$. Moreover, $(G \times
G) \cdot \bm h_J$ is the open dense $G \times G$-orbit in ${\bm X}_J$.

Now for the toroidal embedding $X$ and $K \in \ci$, we may pick a point
$h_K$ in the open $G \times G$-orbit of $X_K$ which maps to  $\bm
h_{p(K)}$. Then $(h_K)_{K \in \ci}$ is a set of representatives of
the $G \times G$-orbits in $X$. Notice that $h_K \in \pi \i(\bm X'_0)
\subset X'_0$.

\subsection{The structure of $G \times G$-orbit closures}

The following result should be well known but, as we have not
been able to find a reference to it, we include a proof.

\begin{lem}
\label{quotient} Let $H$ denote a linear algebraic over the field
$k$ and let $Y$ denote a homogeneous $H$-variety. Let $y \in Y$
and let $p_y : H \rightarrow Y$ denote the associated orbit map.
Let $H_y$ denote the stabilizer group scheme of $y$. Then $Y \simeq
H/H_y$ as homogeneous $H$-varieties. Moreover, if $H_y$ is a
normal subgroup scheme of $H$ then $Y$ may be given a structure of
a linear algebraic group such that $p_y$ is a morphism of
algebraic groups.
\end{lem}
\begin{proof}
By \cite[Prop.III.3.5.2.]{DG} it follows that we may identify
$H/H_y$ with a locally closed subscheme $Y'$ of $Y$. As $Y$ is a
homogeneous $H$-variety we conclude that $Y=Y'$. In particular,
$H/H_y$ is a variety.

Consider now the case when $H_y$ is a normal subgroup scheme of
$H$. Then by \cite[Prop.III.3.5.6]{DG}  the quotient $H/H_y$ is an
affine group scheme. Thus the isomorphism $Y \simeq H/H_y$ induces
a desired algebraic group structure on $Y$.
\end{proof}

\begin{prop}
\label{prop1} Let $X$ be a toroidal embedding. Let $K \in \ci$,
$J=p(K)$ and $h=h_K$. Then
\begin{enumerate}
\item $h$ is invariant under the groups  $U_{I-J}^- \times
U_{I-J}$ and ${\rm diag}(L_{I-J})$. \item We simply write $(L_{I-J}
\times \{1\}) \cdot h$ as $L_{I-J} \cdot h$. The closure
$\overline{L_{I-J} \cdot h}$ in $X$ is $L_{I-J} \times
L_{I-J}$-equivariantly isomorphic to a toroidal equivariant
embedding of a quotient $L_{I-J}/H_{I-J}$ of $L_{I-J}$ by some
(not necessarily reduced) subgroup $H_{I-J}$ of the (scheme
theoretical) center of $L_{I-J}$.
\item The natural morphism
$$ \phi_K : (G \times G) \times_{P_{I-J}^- \times P_{I-J}}
\overline{L_{I-J} \cdot h} \rightarrow X_K,$$
is a birational and bijective $G \times G$-equivariant morphism.
Moreover, when the characteristic of $k$ is positive then $\phi_K$
is an isomorphism.
\item For $v \in W^{I-J}$ and $w \in W$ define
$[K, v, w] := (B \dot{v}, B \dot{w}) \cdot h$.
Then
$$ (G \times G) \cdot h =\bigsqcup_{v \in W^{I-J},
w \in W} [K, v, w].$$
\end{enumerate}
\end{prop}
\begin{proof} The statements holds if $X=\bf X$ (see
\cite[1.1]{Sp}).

Now let $V=\pi^{-1}(\bm h_J)\subset \pi^{-1}({\bf X}'_0 ) = X'_0$. Let
$U_1 \simeq {\mathbb G}_a$ be a 1-dimensional additive subgroup in
$G \times G$ normalized by $T \times T$ which acts trivially on $\bm h_J$. Then
$U_1 \cdot h \subset V \subset X'_0$. By \cite[Prop.6.2.3(ii)]{BK} the
$G \times G$-orbit of $h$ intersects $X'_0$ in a single $T \times T$-orbit.
Hence, $U_1 \cdot h \subset (T \times T) \cdot h$ and thus $U_1$ leaves
$(T \times T) \cdot h$ invariant. But $(T  \times T) \cdot h \simeq (k^*)^n$, for
some $n$, and any action of $ {\mathbb G}_a$ on $(k^*)^n$ is trivial. In
particular, $U_1$ leaves $h$ invariant. This proves that $h$ is
invariant under  $U_{I-J}^- \times U_{I-J}$ and the semisimple
part of ${\rm diag}(L_J)$. Now (1) follows as
any element in the toric variety $X'$ is invariant under
${\rm diag}(T)$.

We identify $X_0$ with $U \times U^- \times {X}'_0$ and simply
write $(T \times \{1\}) \cdot h$ as $T \cdot h$. Then $(U \cap
L_{I-J}) \times (U^- \cap L_{I-J}) \times (\overline{T \cdot h}
\cap X_0')$ is a closed irreducible subset of ${X}_0$ contained in
$\overline{L_{I-J} h}$ and of the same dimension as
$\overline{L_{I-J} h}$. Hence $\overline{L_{I-J} \cdot h} \cap
{X}_0 \simeq (U \cap L_{I-J}) \times (U^- \cap L_{I-J}) \times
(\overline{T \cdot h} \cap X_0')$. As $X'= \overline{T}$ is a
toric variety every $T$-orbit closure in $X'$ is normal. Hence
$\overline{L_{I-J} \cdot  h} \cap {X}_0$ is normal. As a consequence, every
intersection of the form $\overline{L_{I-J} \cdot h} \cap x {X}_0$, for $x
\in L_{I-J}$, is also normal.

We claim that the union $\cup_{x \in L_{I-J}} x X_0$ contains
$\overline{L_{I-J} \cdot h}$. To see this it suffices to prove that the
union $\cup_{x \in L_{I-J}} x {\bf X}_0$ contains the wonderful
compactification $\overline{G_{I-J}} = \overline{L_{I-J} \cdot
\bm h_J}$ (see \cite[1.1]{Sp} for this equality) of $G_{I-J}$.
But ${\bf X}_0$ contains (by definition) the
corresponding open subset $(\overline{G_{I-J}})_0$ of
$\overline{G_{I-J}}$ and, moreover, $\overline{G_{I-J}}$ is
covered by the $L_{I-J}$-translates of the subset
$(\overline{G_{I-J}})_0$. This proves the claim and as a
consequence $\overline{L_{I-J} \cdot h}$ is normal.

As $\pi(h) = \bm h_J$ it follows that the scheme theoretic
$L_{I-J}$-stabilizer $(L_{I-J})_h$ of $h$ is a closed subgroup
scheme of the $L_{I-J}$-stabilizer of $\bm h_J$. The latter stabilizer
coincides with the scheme theoretic center of $L_{I-J}$. So
applying Lemma  \ref{quotient} we conclude that $L_{I-J}
\cdot h$ is isomorphic to the reductive group
$L_{I-J}/(L_{I-J})_h$. As a consequence, $\overline{L_{I-J} \cdot
h}$ is an equivariant embedding of $L_{I-J}/(L_{I-J})_h$. Moreover
the map $\pi$ induces a morphism $\pi : \overline{L_{I-J} \cdot h}
\rightarrow \overline{L_{I-J} \cdot \bm h_J} \simeq
\overline{G_{I-J}}$, so $\overline{L_{I-J} \cdot h}$ is even a
toroidal embedding of $L_{I-J}/(L_{I-J})_h$. This proves statement
(2).

Consider the commutative diagram
$$\xymatrix{
 (G \times G) \times_{P_{I-J}^- \times P_{I-J}} \overline{L_{I-J} \cdot h}
\ar[r]^(0.8){\phi_K} \ar[d] &  X_K \ar[d] \\
(G \times G) \times_{P_{I-J}^- \times P_{I-J}} \overline{L_{I-J}
\cdot \bm h_J}
\ar[r]^(0.8){{\bm \phi}_J} & {\bf X}_J \\
},
$$
where all the maps are the natural ones. As ${\bm \phi}_J$ is an
isomorphism it follows that $\phi_K$ is injective. As $\phi_K$ is
a projective morphism this implies that $\phi_K$ is finite.
Moreover, as $\overline{L_{I-J} \cdot h} $ is closed in $X_K$ and
invariant under ${P_{I-J}^- \times P_{I-J}}$ the image of $\phi_K$
is closed. Therefore $\phi_K$ is surjective and hence bijective.
Moreover, due to the identification  $\overline{L_{I-J} h} \cap
{X}_0 \simeq (U \cap L_{I-J}) \times (U^- \cap L_{I-J}) \times
(\overline{T \cdot h} \cap X_0')$ it follows that $\phi_K$ is
birational. This proves the first part of statement (3). When
the characteristic is positive then $X_K$ is Frobenius split
(see e.g. \cite[Thm.6.2.7]{BK}) and thus weakly normal
(see e.g. \cite[Thm.1.2.5]{BK}). It follows that $\phi_K$
is an isomorphism which ends the proof of statement (3).

By statement (1) and the Bruhat decomposition it easily
follows that the union of $[K, v, w]$, for $v \in W^{I-J},
w \in W$, equals $(G \times G) \cdot h$. Moreover, when
$X= \bm X$ then by \cite[Lemma 1.3(i)]{Sp} this union is
disjoint  (notice that our notation is slightly different
from the notation used in \cite{Sp} : the subset $[J,v,w]$
in \cite{Sp} corresponds to  $[I-J, x, w]$ in the present
paper). As $\pi([K,v,w] )$ equals the associated
$B \times B$-orbit $[J,v,w]$  in $\bm X$ this proves
statement (4) in general.
\end{proof}

\begin{remark}
Statement (3) in Proposition \ref{prop1} above is also
correct in characteristic $0$. This follows from 
Theorem \ref{char0} which proves that $X_K$ is normal 
and thus, by Zariski's main theorem, that $\phi_K$ is 
an isomorphism. A result similar to (3) for some 
special (i.e. {\it regular}) embeddings has earlier 
been obtained in \cite[Sect.2.1]{B}.
\end{remark}

\section{$B \times B$-orbit closures}

In this section we will study inclusions between
$B \times B$-orbit closures in a toroidal embedding $X$ of
$G$. We will state a
precise description of when a $B \times B$-orbit $[K,v,w]$
is contained in the closure of another $B \times B$-orbit
$[K',v',w']$. This generalizes the corresponding
results of T. Springer for $X = \bm X$ given in
\cite[Sect.2]{Sp}. As a consequence we will be able
to prove that any $B \times B$-orbit closure $Z$ in $X$ of
codimension $\geq 2$ is a component of an intersection
of $B \times B$-orbit closures distinct from $Z$.
By standard Frobenius splitting techniques this will
enable us to prove that each $B \times B$-orbit closure
admits a canonical Frobenius splitting.

\subsection{Inclusions between  $B \times B$-orbit closures }
\label{inclusion}
Let $K \in \ci$ and $J=p(K)$. Let $B_J= B \cap L_{I-J}$
and $B'=(w^{I-J}_0 w_0) B (w^{I-J}_0 w_0) \i$. Define $\pi_J: B
\rightarrow B_J$ by $\pi_J(b u)=b$, for $b \in B_J$ and $u \in
U_{I-J}$, and $\pi'_J: B' \rightarrow B_J$ by $\pi'_J (b
u)=b$, for $b \in B_J$ and $u \in U_{I-J}^-$. By Proposition
\ref{prop1}(1) the base point $h_K$ is invariant under ${\rm diag}
(B_J)$. In particular, we may define a $B' \times B$ action on
$G \times G \times \overline{B_J \cdot h_K}$
by
$$(b_1, b_2) (g_1, g_2, z)=(g_1 b_1 \i, g_2 b_2 \i,
(\pi'_J(b_1), \pi_J(b_2)) z),$$
for $b_1 \in B', b_2 \in B$, $g_1, g_2 \in G$ and
$z \in \overline{B_J \cdot h_K}$.
The associated quotient is denoted by
$(G \times G) \times _{B' \times B} \overline{B_J \cdot h_K}$.
The map $G \times G \times \overline{B_J \cdot h_K}
\rightarrow X$, $(g_1, g_2, z) \mapsto (g_1, g_2) z$, induces a
projective surjective morphism 
$$p_K: (G \times G) \times_{B' \times B}
\overline{B_J \cdot h_K} \rightarrow X_K.$$ 
which can be used to prove 

\begin{lem}
\label{incl}
Let $v,v',w,w' \in W$. Assume that $v w_0^{I-J} 
\leq v' w_0^{I-J}$ and $w' \leq w$ in the Bruhat
order on $W$. Then
$$ (B \dot v', B \dot w') \cdot (B_J \cdot h_K) 
\subset \overline{(B \dot v, B \dot w) \cdot  (B_J \cdot h_K)}.$$    
\end{lem}
\begin{proof}
By restricting the map $p_K$ above we obtain a 
projective and surjective  map 
$$(\overline{B \dot v B'} \times \overline{B \dot w B}) 
\times_{B' \times B} \overline{B_J \cdot h_K} \rightarrow 
\overline{(B \dot v, B \dot w) \cdot  (B_J \cdot h_K) }.$$
For the above statement to be true
it thus suffices to have $B v' B' \subset \overline{B v B'}$ 
and $B w' B \subset \overline{B w B}$, which is clearly
satisfied under the stated conditions. 
\end{proof}

Notice that when $v \in W^{I-J}$ then the set 
$(B \dot v, B \dot w)  \cdot (B_J \cdot h_K)$,
in Lemma \ref{incl}, coincides with the orbit 
$[K,v,w]$.

\begin{prop}
\label{Bclosure}
Let $K, K' \in \ci$, $v \in W^{I-p(K)}$,
$v' \in W^{I-p(K')}$ and $w, w' \in W$. Then $[K', v',
w']$ is contained in $\overline{[K, v, w]}$ if and only if $K
\subset K'$ and there exists $u \in W_{I-p(K')}$ and $u' \in
W_{I-p(K)} \cap W^{I-p(K')}$ such that $v u' u \i \le v'$, $w' u
\le w u'$.
\end{prop}
\begin{proof} Notice $\overline{[K, v, w]} \subset \pi
\i(\pi(\overline{[K, v, w]})) \cap X_K$. Thus if $[K',
v', w'] \subset \overline{[K, v, w]}$, then $K \subset K'$ and
$[p(K'), v', w'] \subset \overline{[p(K), v, w]}$. By
\cite[2.4]{Sp}, there exists $u \in W_{I-p(K')}$ and $u' \in W_{I-p(K)}
\cap W^{I-p(K')}$ such that $v u' u \i \le v'$, $w' u \le w u'$.

On the other hand, assume that $v' \in W^{I-p(K')}$, $w' \in W$, $u
\in W_{I-p(K')}$ and $u' \in W_{I-p(K)} \cap W^{I-p(K')}$ such that
$v u' u \i \le v'$, $w' u \le w u'$. Assume, moreover, that $K
\subset K'$. By the one to one correspondence between the set
of $G \times G$-orbits in $X$ and the set of $T$-orbits
in $X_0'$ \cite[Prop.6.2.3(ii)]{BK}, it follows that
$h_{K'} \in \overline{T \cdot h_K}$. Thus $(\dot x, \dot x) h_{K'} \in
\overline{T \cdot h_K}$ for all $x \in W_{I-p(K)}$ by Proposition
\ref{prop1}(i). Therefore, with $J' = p(K')$, we find by use
of Lemma \ref{incl},
$$[K', v', w']=(B \dot v' \dot u, B \dot w' \dot u) \cdot h_{K'}
\subset \overline{(B \dot v \dot u', B \dot w \dot u') \cdot 
(B_{J'} \cdot h_{K'})}.$$
As $u' \in W^{I-J'}$ we have $u' B_{J'} \subset B u'$.  
Thus the right hand side of the above inclusion is contained in 
$$ \overline{(B \dot v B , B \dot w ) \cdot ( (u',u') h_{K'})}  
\subset \overline{(B \dot v , B \dot w ) \cdot(B_J \cdot h_{K})}
= \overline{[K,v,w]},$$  
which ends the proof.
\end{proof}

We may reformulate the above proposition to a slightly simpler
version.

\begin{prop}
\label{reformulation}
Let $K, K' \in \ci$, $v \in W^{I-p(K)}$, $v' \in W^{I-p(K')}$ and
$w, w' \in W$. Then $[K', v', w'] \subset \overline{[K, v, w]}$ if
and only if $K' \supset K$ and there exists $u \in W_{I-p(K)}$
such that $v u \le v'$, $w' \le w u$.
\end{prop}
\begin{proof} If $[K', v', w'] \subset \overline{[K, v, w]}$, then in
$\bm X$ we have $$[I, v', w'] \subset \overline{[p(K'), v', w']}
\subset \overline{[p(K), v, w]}.$$ 
By Proposition \ref{Bclosure}
there exists $u \in W_{I-p(K)}$ such that $v u \le v'$, $w' \le w
u$. On the other hand, assume that $K' \supset K$ and there exists 
$u \in W_{I-p(K)}$ such that $v u \le v'$, $w' \le w u$. 
Write $u$ as $u=u_1 u_2$
for $u_1 \in W_{I-p(K)} \cap W^{I-p(K')}$ and $u_2 \in
W_{I-p(K')}$. By \cite[Cor.3.4]{He}, there exists $u'_2 \leq u_2$
such that $w' (u'_2) \i \le w u_1$. Moreover, $v u_1 u'_2 \leq v
u_1 u_2 \leq v'$. Hence by Proposition \ref{Bclosure}, $[K', v',
w'] \subset \overline{[K, v, w]}$ and the proposition is proved.
\end{proof}

For later reference we state the following easy consequences
of the above propositions.

\begin{cor}
\label{consequence}
Let $K, K' \in \ci$, $v \in W^{I-p(K)}$, $v' \in W^{I-p(K')}$ and
$w, w' \in W$.
\begin{enumerate}
\item If  $\overline{[K', v', w']} \subset \overline{[K, v,
w]}$ then $v \leq v'$.
\item $\overline{[K, v, w']} \subset \overline{[K', v, w]}$
if and only if  $w' \leq w$ and $K' \subset K$.
\end{enumerate}
\end{cor}

\subsection{Intersection of $B \times B$-orbit closures}

In this section we will prove.

\begin{prop}
\label{componentintersection}
Let $Z \neq X$ denote a $B \times B$-orbit closure in $X$.
If $Z$ has codimension 1 in $X$ the $Z$ is either a boundary
divisor $X_i$, $1 \le i \le n$, of $X$ or else $Z$ coincides
with the closure of a
codimension 1 Bruhat cell $B \dot{s_i}\dot{w_0} B$, $1 \le i
\le l$,  within
$X$. If the codimension of $Z$ is $\ge 2$ then there exist
$B \times B$-orbit closures $Z_1 \neq Z$ and $Z_2 \neq Z$
in X such that $Z$ is a component of the intersection
$Z_1 \cap Z_2$.
\end{prop}

The proof of Proposition \ref{componentintersection} will
depend on the following 4 lemmas.

\begin{lem}
\label{comp1} Let $w \in W$ be an element of length $l(w) < l(w_0)
-1$. Then there exist elements $w'$ and $w''$ distinct from $w$
such that $\overline{[\varnothing, 1, w]}$ is an irreducible
component of $\overline{[\varnothing, 1, w']} \cap
\overline{[\varnothing, 1, w'']}$.
\end{lem}
\begin{proof}
Choose simple reflections $s_i$ and $s_j$ such that $l(w s_i) =
l(s_j w) = l(w)+1$. If $w s_i$ and $s_j w$ are distinct then the
statement follows by setting $w' = w s_i$ and  $w'' = s_j w$. If
$w s_i = s_j w$, then we choose a simple reflection $s_k$ such
that $l(w s_i s_k)=l(w s_i)+1=l(w)+2$. Then $k \neq i$. As $w s_i
s_k=s_j w s_k$, we conclude that $l(w s_k)=l(w)+1$. The statement
follows by setting $w'=w s_i$ and $w''=w s_k$.
\end{proof}

\begin{lem}
\label{comp2} For $K \in \ci$ and $w \in W$, $\overline{[K, 1,
w]}$ is an irreducible component of $\overline{[\varnothing, 1,
w]} \cap \overline{[K,1,w_0]} $.
\end{lem}
\begin{proof} By Proposition \ref{reformulation},
$\overline{[K, 1, w]} \subset \overline{[\varnothing, 1,
w]} \cap \overline{[K,1,w_0]}$. As $X$ is a finite union
of $B \times B$-orbits each irreducible component of
the intersection $\overline{[\varnothing, 1, w]}
\cap  \overline{[K, 1,w]} $ will be the closure of a
$B \times B$-orbit in $X$.
Assume that $K' \in \ci$, $v \in W^{I-p(K')}$ and $w' \in
W$ satisfy
$$\overline{[K, 1, w]} \subset \overline{[K', v, w']}
\subset \overline{[\varnothing, 1, w]}
\cap \overline{[K,1,w_0]}.$$
Then by Corollary \ref{consequence}(1) we have $v =1$.
Moreover, Proposition \ref{reformulation} implies that
$K'=K$. Then Corollary \ref{consequence}(2) shows that $w'=w$,
which ends the proof.
\end{proof}

\begin{lem}
\label{comp3} Let $v, v' \in W^{I-p(K)}$ with $v=s_i v'$ for some
$i \in I$ and $l(v)=l(v')+1$. Then $\overline{[K, v, w_0]}$ is an
irreducible component of $\overline{[K, v', w_0]} \cap
\overline{[\varnothing, 1, w_0 v \i]}$.
\end{lem}
\begin{proof}
By Proposition  \ref{reformulation} we easily conclude
 $\overline{[K, v, w_0]} \subset
\overline{[\varnothing, 1, w_0 v \i]}$ and $\overline{[K, v, w_0]}
\subset \overline{[K, v', w_0]}$. Assume that $w \in W^{I-p(K)}$
and $w' \in W$ satisfy
$$\overline{[K, v, w_0]} \subset \overline{[K, w, w']}
\subset \overline{[K, v', w_0]} \cap \overline{[\varnothing, 1, w_0 v \i]}.$$
Then by Corollary
\ref{consequence} (i), $v' \le w \le v$. So $w=v'$ or $w=v$.
Moreover, by Proposition \ref{reformulation} there exists
$u \in W_{I-p(K)}$ such that $w u  \le v$ and $w_0  \le w' u$.
As $v \in W^{I-p(K)}$ we
conclude that $u=1$ and $w'= w_0$. Then, by Proposition
\ref{reformulation}, there
exists $u' \in W$ such that $u' \leq w$ and $w_0 \leq w_0 v \i
u'$. Thus $u' = v$ and $w$ must then be equal to $v$. The lemma is
proved.
\end{proof}

\begin{lem}
\label{comp4} We keep the assumptions on $v$ and $v'$ from the
previous Lemma \ref{comp3}. Then for $w \in W$, $\overline{[K, v, w]}$
is an irreducible component of $\overline{[K, v, w_0]} \cap
\overline{[K, v', w]}$.
\end{lem}
\begin{proof}
By Proposition \ref{reformulation} we have
 $\overline{[K, v, w]} \subset \overline{[K, v,
w_0]} \cap \overline{[K, v', w]}$. Assume that $u \in W^{I-p(K)}$
and $w' \in W$ satisfy $\overline{[K, v, w]} \subset \overline{[K,
u, w']} \subset \overline{[K, v, w_0]} \cap \overline{[K, v', w]}$.
Then, by Corollary \ref{consequence}(i), $u=v$ and hence by
Corollary \ref{consequence}(ii) we have $w \leq w'$. Moreover, by
Proposition \ref{reformulation} there exists $u' \in W_{I - p(K)}$ such
that $v' u' \leq v$ and
$w' \leq w u'$. We conclude that $u=1$ and as a consequence that $w'
=w$.
\end{proof}

We can now prove  Proposition \ref{componentintersection}.

\begin{proof}
Let $K \in \ci$, $v \in W^{I-p(K)}$ and $w \in W$ such that
$Z = \overline{[K,v,w]}$. Notice that by Proposition
\ref{reformulation} the closure $\overline{[\varnothing,1,w_0]}$
contains all $B \times B$-orbit closures and hence it will
be equal to $X$.

We first consider the situation when $w \neq w_0$ :
if there exists a simple reflection $s_i$
such that $l(s_i v) = l(v)-1$ then by Lemma \ref{comp4}
we may use $Z_1 =  \overline{[K, v, w_0]}$ and $Z_2 =
\overline{[K, s_i v, w]}$ (notice that this makes sense
as $s_i v \in W^{I-p(K)}$). So we may assume that $v=1$.
If now $K \neq \varnothing$ then by Lemma \ref{comp2}
we may use $Z_1 = \overline{[\varnothing, 1, w]}$ and
$ Z_2 =\overline{[K,1,w_0]}$. So we may assume that 
$Z=\overline{[\varnothing, 1, w]}$. If $l(w) <
l(w_0)-1$ then we may apply Lemma \ref{comp1} to define
$Z_1$ and $Z_2$. This leaves us with the cases
$\overline{[\varnothing, 1, s_i w_0 ]}$, $i=1, \dots,l$,
which are equal to the closures of the Bruhat cells
$B \dot{s_i}\dot{w_0} B \subseteq G$ within $X$.

Next assume that $w = w_0$ : if there exists a simple
reflection $s_i$ such that $l(s_i v) = l(v)-1$ then by
Lemma \ref{comp3} we may use $Z_1 = \overline{[K, s_iv, w_0]}$
and $Z_2 = \overline{[\varnothing, 1, w_0 v \i]}$. So
we may assume that $v=1$. As $Z \neq X$ we have
that $Z = \overline{[K,1,w_0]}$ with $K$ a nonempty set.
In particular, by Proposition \ref{reformulation}, $Z$
coincides with the $G \times G$-orbit closure $X_K$. 
Now let $K' \subset K$ be a minimal subset such that 
$X_K$ is an irreducible component of $\cap_{i \in K'} X_i$.
If $|K'| =1 $ then $Z$ coincides with a boundary 
divisor, so we may assume that $|K'| > 1 $.
Let $j \in K'$ and let $Y_1, \dots, Y_s$ denote the 
irreducible components of the intersection $\cap_{i 
\in K'-\{ j \} } X_i$. Then, by minimality of $K'$, 
each $Y_i$, for $i=1, \dots,s $, is a $G \times G$-orbit
closure distinct from $Z$. Moreover, there exists an
$i$ such that $Z$ is an irreducible 
component of the intersection $Y_i \cap X_j$. Now 
use $Z_1 = Y_i$ and $Z_2 = X_j$. 
\end{proof}

\section{Frobenius splitting of $B \times B$-orbit
closures}

Let $X$ denote an equivariant embedding of the reductive
group $G$ over a field of positive characteristic $p>0$.
As above the boundary divisors of $X$ will be denoted by
$X_1, \dots, X_n$. Moreover, we will use the notation
$D_i$, $i =1, \dots,l$, to denote the closures of the
Bruhat cells $B \dot{s_i} \dot{w_0} B$, $i=1, \dots,l$, within $X$.

\begin{prop}
\label{Frobeniussplitting}
The equivariant embedding $X$ admits a $(B \times B, T
\times T)$-canonical Frobenius splitting which compatibly
splits the closure of all $B \times B$-orbits.
\end{prop}
\begin{proof}
First of all $X$ admits a $(B \times B, T \times T)$-canonical
Frobenius splitting $s$ which compatibly splits all boundary
component $X_j$, $j=1, \dots, n$, and the subvarieties
$D_i$, $i=1, \dots, l$ (see \cite[Thm.6.2.7]{BK}).

Consider, for a moment, the case when $X$ is toroidal. We claim
that $s$ compatibly Frobenius splits all $B \times B$-orbit
closures. If this is not the case, then there exists a
$B \times B$-orbit closure $Z$ of maximal dimension which is
not compatibly Frobenius split by $s$. By Proposition
\ref{componentintersection} the codimension of $Z$ must
be $\ge 2$. In particular, we can find orbit closures
$Z_1 \neq Z$ and $Z_2 \neq Z$ such that $Z$ is a component
of the intersection $Z_1 \cap Z_2$. By the maximality assumption
on $Z$ the orbit closures $Z_1$ and $Z_2$ will be compatibly
Frobenius split by $s$. But then every component of $Z_1 \cap Z_2$, 
and thus $Z$, will also be compatibly Frobenius
split by $s$, which is a contradiction. This ends the proof
when $X$ is toroidal.

For an arbitrary embedding $X$ we may find a toroidal embedding
$X'$ of $G$ and a birational projective morphism
$f : X' \rightarrow X$  extending the identity map on $G$
(see e.g. \cite[Prop.6.2.5]{BK}). Now $X'$ admits a
$(B \times B, T \times T)$-canonical Frobenius splitting
$s'$ which compatibly Frobenius splits all $B \times B$-orbit
closures. By Zariski's main theorem the map
$f^\sharp : \mathcal{O}_{X'} \rightarrow f_* \mathcal{O}_{X}$
induced by $f$ is an isomorphism. In particular, $s'$ induces
by push forward a $(B \times B, T \times T)$-canonical Frobenius
splitting $s$ of $X$. Moreover, the image in $X$ of every
$B \times B$-orbit closure in $X'$ will be compatibly Frobenius
split by $s$. But any $B \times B$-orbit closure in $X$
is the image of a similar orbit closure in $X'$. This ends
the proof.
\end{proof}

\subsection{Cohomology vanishing}

As a direct consequence of Proposition
\ref{Frobeniussplitting} we conclude
the following vanishing result (see
e.g. \cite[Thm.1.2.8]{BK}).

\begin{prop}
\label{vanishing}
Let $X$ be a projective equivariant embedding of $G$.
Let $Z$ denote a $B \times B$-orbit closures in $X$
and let $\mathcal L$ denote an ample line bundle on
$Z$. Then
$$ {\rm H}^i\big(Z, \mathcal L \big)=0 , ~ i>0.$$
Moreover, if $Z' \subset Z$ is another $B \times B$-orbit
closure then the restriction map
$${\rm H}^0\big( Z, \mathcal L\big)
\rightarrow  {\rm H}^0\big(Z', \mathcal
L_{Z'}\big),$$
is surjective.
\end{prop}

Later (Corollary \ref{nef}) we will see that the 
vanishing part of Propositione \ref{vanishing}
remains true when the line bundle $\mathcal L$ 
is only assumed to be nef, i.e. when $\mathcal L 
\otimes \mathcal M$ is an ample line bundle for 
every ample line bundle $\mathcal M$.

\section{Global $F$-regularity of $B \times B$-orbit closures}

We are now ready to state and prove the main result of the
paper.

\begin{thm}
\label{glreg}
Let $X$ denote a projective equivariant embedding
of a reductive group $G$ over a field of positive
characteristic $p>0$. Let $Z$ denote a $B \times B$-orbit
closure in $X$. Then $Z$ is globally $F$-regular.
\end{thm}

We will divide the proof of Theorem \ref{glreg} into
2 parts. The first part concerns the case when $X$
is toroidal.

\begin{lem}
\label{glreg1}
Let $X$ be a projective toroidal embedding.
Then any $B \times B$-orbit closure $\overline{[K, v, w]}$
in $X$ is globally $F$-regular.
\end{lem}
\begin{proof}
Keep the notation of Section \ref{inclusion}. As a consequence
of Proposition \ref{Frobeniussplitting}, $X$ admits a
$(B' \times B, T \times T)$-canonical Frobenius splitting $s$
which compatibly Frobenius splits every $B' \times B$-orbit
closure.

Let $Y = \overline{(B' \times B) h_K}$ and  $Y' = Y -
(B' \times B) h_K$. Then  $s$ induces a
$(B' \times B, T \times T)$-canonical Frobenius splitting
$s_Y$ of $Y$ which compatibly Frobenius splits $Y'$.
Notice that by Proposition \ref{prop1}(1),
$ Y = \overline{B_J \cdot h_K}$. Thus by Proposition
\ref{prop1}(2), $Y$ is
the closure of the Borel subgroup $B_J$ of $L_{I-J}$ within
some equivariant embedding of $L_{I-J}$. Hence, $Y$ is a large
Schubert variety for some equivariant embedding of $L_{I-J}$
and, as such, $Y$ is globally $F$-regular \cite[Thm.4.3]{BT}.
Define $v'=w^{I-p(K)}_0 w_0 v$. Then $(B' \times B)_{(v',w)}$
contains the group $B_J \times \{ 1 \}$ (notice that the set
of positive roots on the first coordinate is defined with
respect to $B'$) and thus by  Proposition \ref{prop1}(1)
$$(B' \times B) h_K = (B_J \times \{ 1 \}) h_K =
(B' \times B)_{(v',w)} h_K.$$

The above statements proves that the triple
$(Y, h_K , (v', w))$ satisfies the requirements of
Proposition \ref{mainlemma}. Now Theorem \ref{glregthm}
shows that the closed subvariety
$$(\overline{B' \dot v' B'}, \overline{B \dot w B}) \overline{B_J \cdot
h_K}=(\dot w^{I-J}_0 \dot w_0 , 1) \overline{[K, v, w]},$$
is globally $F$-regular. Thus also $ \overline{[K, v, w]}$
must be globally $F$-regular.
\end{proof}

\subsection{The general case}

Let $X$ denote an arbitrary equivariant projective
embedding of $G$. To handle the proof of Theorem 
\ref{glreg}
for $X$ we start by the following construction :
Consider the natural $G \times G$-equivariant embedding
$$  f : G \rightarrow X \times {\bf X}.$$
and let $Y$ denote the normalization of the closure
of the image of $f$. Then $Y$ is a projective
equivariant toroidal embedding of $G$. We let
$\phi : Y \rightarrow X$ denote the associated
$G \times G$-equivariant projective morphism to $X$.
Then

\begin{lem}
\label{glreg2} Let $Z'$ denote the closure of a $B \times
B$-orbit within $Y$ and let $Z$ denote its image $\phi(Z')$
within $X$. Then the induced morphism $\phi' : Z' \rightarrow
Z$ is a rational morphism.
\end{lem}
\begin{proof}
We will prove this using Lemma \ref{Kempf}.
Notice first of all that $\phi$ is birational
and $X$ is normal, so by Zariski's main theorem
we have $\phi_* \mathcal O_{Y} = \mathcal O_X$.
Let now $\mathcal L$ denote a very ample line bundle on $X$.
Then by Lemma \ref{glreg1} and \cite[Cor.4.3]{S},
$$ {\rm H}^i\big(Y, \phi^* \mathcal L\big) =
{\rm H}^i\big(Z', \phi^* \mathcal L\big) = 0, ~i>0,$$
as $ \phi^* \mathcal L$ is globally generated
and thus nef.

Let $\tilde{\bm D}_i$, $i=1, \dots, l$, denote the closures
$\overline{B^- \dot{s}_i \dot{w}_0 B^-}$ in $\bm{X}$. Then the divisor
$\tilde{\bm D} = \sum_{i=1}^l  \tilde{\bm D}_i$ is ample
\cite[Prop.6.1.11]{BK}. Let $\mathcal M = \mathcal
O_{\bm X}(\tilde{\bm D})$ denote the associated line bundle
and let $\mathcal M' = \phi^* \mathcal M$ be its pull back
to $Y$. Let $\bm s$ denote the canonical section of
$\mathcal M$ and let $s'$ denote its pull back to $Y$.
Let $V$ denote an irreducible component of the
support of $s'$. If $V$ is contained in the boundary of
$Y$ then the support of $s'$ will contain a closed
$G \times G$-orbit. In particular, also the support
$\cup_i \tilde{\bf D}_i$ of $s$ will contain a closed
$G \times G$-orbit. As the latter is not the case
we conclude that each component of the support of $s'$ will
intersect $G$. Moreover, the support of $s'$ is $B^- \times
B^-$-stable. As a consequence, we conclude that the divisor of
zeroes of $s'$ equals
$$\sum_{i=1}^l n_i \tilde{D}'_i,$$
for some positive integers $n_i$ and with $\tilde{D}'_i$, $i=1,
\dots, l$, denoting the closure $\overline{B^- \dot s_i \dot w_o B^-}$ in $Y$. 

Let $Y_j$, $j=1, \dots, n$, denote the boundary components in
$Y$ and let $D_i'$, $i=1, \dots, l$, denote the closures
$\overline{B \dot{s}_i \dot{w}_0 B}$ in $Y$. Let $Y^0$
denote the smooth locus of $Y$. Then $Y^0$ admits a
Frobenius splitting which compatibly Frobenius splits the
Cartier divisors $Y^0 \cap Y_j$, $j=1, \dots, n$, and
$D_i' \cap Y^0$ and $\tilde{D}_i' \cap Y^0$, $i=1,
\dots,l$ \cite[Thm.6.2.7]{BK}. As in the proof of Proposition \ref{rational}
we conclude that $Y^0$ admits a stable Frobenius
splitting along the effective divisor
$$ {\rm div}(s') \cap Y^0 = \sum_{i=1}^l n_i (\tilde{D}'_i \cap Y^0),$$
which compatibly Frobenius splits $D_i' \cap Y^0$, $i=1,
\dots,l$, and $Y_j \cap Y^0$, $j=1, \dots, n$. Let $\psi_0$
denote such a stable Frobenius splitting; i.e. let $e$ be an
integer such that $\psi_0$ is a splitting of the morphism
$$ \mathcal O_{Y^0} \rightarrow  F_*^e\mathcal
M'_{| Y^0},$$ defined by the restriction of $s'$ to $Y^0$.
Let now $i : Y^0 \rightarrow Y$ denote the inclusion morphism.
Applying the functor $i_*$ to the above split morphism and using
that $Y$ is normal, we find that the morphism
$$  \mathcal O_{Y} \rightarrow  F_*^e \mathcal
M',$$ defined by $s'$ has an induced  splitting $\psi$. Then
$\psi$ defines a stable Frobenius splitting along ${\rm div}(s')$
which compatibly Frobenius splits $D_i'$, $i=1, \dots,l$, and
$Y_j$, $j=1, \dots, n$ (as the compatibility can be checked on
the open dense subsets $Y^0$).

We now claim that $Z'$ is not contained in any $\tilde{D}_i'$. To
see this assume that $Z'$ is contained in  $\tilde{D}_i'$ for some
$i$. As $Z'$ is $B \times B$-invariant and as $\tilde{D}_i'$ is
$B^- \times B^-$-invariant it follows that $(B^-B, B^-B) Z'$ is
contained in $\tilde{D}_i'$. But then also $(G,G)Z'$ must be
contained in $\tilde{D}_i'$. We conclude that
$\tilde{\bf D}_i'$ contains a closed $G \times G$-orbit
which is a contradiction. Hence, $Z'$ is not contained in the
support of $s'$. As in the proof of Proposition
\ref{Frobeniussplitting} we may then use Proposition
\ref{componentintersection} to show that $Z'$ is
compatibly Frobenius split by the stable Frobenius
splitting $\psi$. By \cite[Lem.4.8]{T} it follows
that we have an embedding
$${\rm H}^1\big(Y, \mathcal I_{Z'} \otimes  \phi^* \mathcal
L \big) \subset
{\rm H}^1\big(Y, \mathcal I_{Z'} \otimes  \phi^* \mathcal
L^{p^e} \otimes \mathcal M' \big) $$
of abelian groups, where $\mathcal I_{Z'}$ denotes
the sheaf of ideals associated to $Z'$. But
$\mathcal L^{p^e} \boxtimes \mathcal M$ is ample
on $X \times \bm X$ and, as the map $Y \rightarrow
X \times \bm X$ is finite, we conclude that
$ \phi^* \mathcal L^{p^e} \otimes \mathcal M'$
is ample on $Y$. Applying  \cite[Thm.1.2.8]{BK}
it follows that $ {\rm H}^1\big(Y, \mathcal I_{Z'}
\otimes  \phi^* \mathcal L \big)$ is zero.

As all the requirement in Lemma \ref{Kempf} are now
satisfied this ends the proof.
\end{proof}

We may now prove  Theorem \ref{glreg}

\begin{proof}
By Corollary \ref{pushforward} and Lemma \ref{glreg2} we may
assume that $X$ is toroidal. Now apply Lemma \ref{glreg1}.
\end{proof}

\subsection{Applications}

As the main application of  Theorem \ref{glreg} we
find.

\begin{cor}
\label{glregcor}
Let $X$ denote an equivariant embedding of a reductive
group $G$ over a field of positive characteristic. Then
every $B \times B$-orbit closure in $X$ is strongly
$F$-regular. In particular, every $B \times B$-orbit
closure is normal, Cohen-Macaulay and locally
$F$-rational.
\end{cor}
\begin{proof}
As in the proof of \cite[Cor.4.2]{BT} we may reduce
to the case when $X$ is projective. Then by Theorem
\ref{glreg} every $B \times B$-orbit closure is
globally $F$-regular and thus strongly $F$-regular.
This ends the proof.
\end{proof}

We also obtain the following strengthening of Proposition
\ref{vanishing}.

\begin{cor}
\label{nef}
Let $X$ denote a projective equivariant embedding of a
reductive group $G$ over a field of positive characteristic.
Let $Z$ denote the closure of a $B \times B$-orbit and let
$\mathcal L$ be a nef line bundle on $Z$. Then the cohomology
${\rm H}^i \big ( Z, \mathcal L \big)$ vanishes for $i>0$.
\end{cor}
\begin{proof}
Just apply \cite[Cor.4.3]{S}.
\end{proof}

\section{The characteristic 0 case}

Let $X$ denote a scheme of finite type over a field $K$ of
characteristic $0$. Then there exists a finitely generated
$\mathbb Z$-algebra $A$ and a flat scheme $X_A$ of finite type
over $A$, such that the base change of $X_A$ to $K$ may be
naturally identified with $X$. Moreover, when $m \subset A$ is a
maximal ideal we may form the base change $X_{k(m)}$ of $X_A$ to
the finite field $k(m)=\nicefrac{A}{m}$. We then say that the
scheme $X$ is of {\it strongly $F$-regular type} (resp. {\it
$F$-rational type}) if $X_{k(m)}$ is strongly $F$-regular (resp.
$F$-rational) for all maximal ideals $m$ in a dense open subset of
${\rm Spec(A)}$.

Any scheme $X$ of strongly $F$-regular type will also be of
$F$-rational type. Thus, by \cite[Thm.4.3]{S2}, schemes of strongly
$F$-regular type will have rational singularities, in particular,
they will be normal and Cohen-Macaulay.

In the proof of the next result we will use the following
observation (see e.g. \cite[Thm.5.5(e)]{HH3}: let
$\overline{k}(m)$ denote an algebraic closure of the field $k(m)$.
If the base change $X_{\overline{k}(m)}$ is strongly $F$-regular
then also $X_{k(m)}$ is strongly $F$-regular.

We can now prove the characteristic $0$ version of Corollary
\ref{glregcor}.

\begin{thm}
\label{char0}
Let $X$ denote an equivariant embedding of a reductive group $G$
over an algebraically closed field $k$ of characteristic $0$. Then
every $B \times B$-orbit closure in $X$ is of strongly $F$-regular
type. In particular, every $B \times B$-orbit closure in $X$ has
rational singularities.
\end{thm}
\begin{proof}
We may assume that there exists a split $\mathbb{Z}$-form
$G_{\mathbb Z}$ of $G$ over which $B$ is defined by a closed
subscheme $B_{\mathbb Z}$. Let $Z$ denote a $B \times B$-orbit
closure in $X$. The complete data consisting of the $G \times
G$-action on $X$, the open embedding $G \subset X$, the $B \times
B$-stability of $Z$, the closed embedding $Z \subset X$ and the
irreducibility of $X$ and $Z$ may all be descended to some
finitely generated $\mathbb Z$-algebra $A$ (see e.g.
\cite[Sect.2]{HH2} for this kind of technique). This means that
there exists schemes $G_A := G_{\mathbb Z} \times_{{\rm
Spec}(\mathbb Z)} {\rm Spec}(A)$, $B_A := B_{\mathbb Z}
\times_{{\rm Spec}(\mathbb Z)} {\rm Spec}(A)$, $X_A$ and $Z_A$
flat and of finite type over ${\rm Spec}(A)$ satisfying, that for
every maximal ideal $m \subseteq A$ the associated base changes
$G_{\overline{k}(m)}$,
 $B_{\overline{k}(m)}$, $X_{\overline{k}(m)}$ and
$Z_{\overline{k}(m)}$, to an algebraic closure $\overline{k}(m)$
of the field $k(m)=\nicefrac{A}{m}$, share the same structure;
i.e. $G_{\overline{k}(m)}$ is a reductive linear algebraic group,
$X_{\overline{k}(m)}$ is an irreducible $G_{\overline{k}(m)}
\times G_{\overline{k}(m)}$-variety containing
$G_{\overline{k}(m)}$ as an open subset and $Z_{\overline{k}(m)}$
is an irreducible $B_{\overline{k}(m)} \times
B_{\overline{k}(m)}$-stable subvariety of $X_{\overline{k}(m)}$.
As $X$ is normal we may even assume that $X_{\overline{k}(m)}$ is
normal  (see \cite[Thm.2.3.17]{HH2}). In particular,
$X_{\overline{k}(m)}$ is then an equivariant embedding of the
reductive group $G_{\overline{k}(m)}$. Moreover, by the finiteness
of the number of $B_{\overline{k}(m)} \times
B_{\overline{k}(m)}$-orbits in $X_{\overline{k}(m)}$ we conclude
that  $Z_{\overline{k}(m)}$ is the closure of such an orbit.

Applying Corollary \ref{glregcor} and the observation above, 
we conclude that $Z$ is of strongly $F$-regular type and thus
also of $F$-rational type. Finally, as mentioned above, the 
latter statement implies that $Z$ has rational singularities.
\end{proof}

We may now generalize Corollary \ref{nef} to arbitrary 
characteristics.

\begin{cor}
Let $X$ denote a projective equivariant embedding of a
reductive group $G$ over a field of arbitrary characteristic. 
Let $Z$ denote the closure of a 
$B \times B$-orbit and let $\mathcal L$ be a nef line 
bundle on $Z$. Then the cohomology ${\rm H}^i \big ( 
Z, \mathcal L \big)$ vanishes for $i>0$.
\end{cor}
\begin{proof}
Apply Corollary \ref{nef} and \cite[Cor.5.5]{S}.
\end{proof}
  
For a discussion of other kinds of vanishing results 
for varieties of globally $F$-regular type we refer 
to \cite[Sect.5]{S}.


\begin{thebibliography}{WWWW}
  \bibitem[B-K]{BK} M. Brion and S. Kumar, \emph{Frobenius
    Splittings Methods in Geometry and Representation Theory},
    Progress in Mathematics (2004), Birkh\"auser, Boston.

\bibitem[B]{B} M. Brion, \emph{The behaviour at infinity of the
Bruhat decomposition}, Comment. Math. Helv. \textbf{73}  (1998),
137--174.

\bibitem[B2]{B2} M. Brion, \emph{Multiplicity-free subvarieties 
of flag varieties}, in {\it Commutative Algebra}, 
Contemporary Math. \textbf{331} (2003) 13--23.


  \bibitem[B-P]{BP} M. Brion and P. Polo, \emph{Large Schubert
    Varieties}, Represent. Theory  \textbf{4}  (2000), 97--126.

  \bibitem[B-T]{BT} M. Brion and J. F. Thomsen, \emph{F-regularity of
      large Schubert varieties}, math.AG/0408180.


  \bibitem[D-G]{DG} M. Demazure and P. Gabriel, \emph{Groupes alg\'ebriques. {T}ome {I}:
    {G}\'eom\'etrie alg\'ebrique, g\'en\'eralit\'es, groupes
    commutatifs}, Masson \& Cie, 1970.

\bibitem[H]{Hartshorne} R~.Hartshorne, \emph{Ample Subvarieties of
    Algebraic Varieties}, Lecture notes in math. \textbf{156} (1970),
    Springer-Verlag.

  \bibitem[H-H]{HH} M.~Hochster and C.~Huneke, \emph{{Tight closure
        and strong $F$-regularity}},
    M\'emoires de la Soci\'et\'e Math\'ematique de France \textbf{38} (1989),
    119--133.

  \bibitem[H-H2]{HH2} M.~Hochster and C.~Huneke, \emph{F-regularity,
test elements, and smooth base change}, Trans. Amer. Math. Soc.
\textbf{346} (1994), 1--62.

 \bibitem[H-H3]{HH3} M.~Hochster and C.~Huneke, \emph{Tight
closure in equal characteristic zero}, preprint.


  \bibitem[He]{He} X. He, \emph{The $G$-stable pieces of the wonderful
      compactification}, math.RT/0412302.

\bibitem[L-S]{LS}
    G.~Lyubeznik and K.~E.~Smith, \emph{Weak and strong $F$-regularity
      are equivalent in graded rings}, Amer. Math. Soc.
    \textbf{121} (1999), 1279--1290.

 \bibitem[L-T]{LT} N.~Lauritzen and J.~F.~Thomsen, \emph{Line
bundles on Bott-Samelson varieties}, J. Alg. Geom.
\textbf{13} (2004), 461--473.


  \bibitem[L-P-T]{LPT} N.~Lauritzen, U.~R.~Pedersen and J.~F.~Thomsen, \emph{Global
     $F$-regularity of Schubert varieties with applications to
     $D$-modules}, math.AG/0402052. 


 \bibitem[R]{R} A,~Ramanathan, \emph{Equations defining Schubert varieties and Frobenius
      splitting of diagonals}, Inst. Hautes \'Etudes Sci. Publ. Math. \textbf{65} (1987), 61--90.

 \bibitem[S]{S2} K.~E.~Smith, \emph{$F$-rational rings have rational
       singularities}, Amer. J. Math. \textbf{119} (1997), 159--180.

 \bibitem[S2]{S}
    K.~E.~Smith, \emph{{Globally $F$-regular varieties: Applications
    to vanishing theorems for quotients of Fano varieties}}, Michigan
    Math.~J. \textbf{48} (2000), 553--572.

 \bibitem[Sp]{Sp} T.~A.~Springer, \emph{Intersection cohomology of $B \times
    B$-orbits closures in group compactifications}, J. Alg. {\bf 258}
    (2002), 71--111.

\bibitem[T]{T} J.F. Thomsen, \emph{Frobenius splitting of
      equivariant closures of regular conjugacy classes},
    math.AG/0502114.

\end{thebibliography}
\end{document}